\documentclass[journal]{IEEEtran}

\usepackage{graphics} % for pdf, bitmapped graphics files
\usepackage{epsfig} % for postscript graphics files
\usepackage{times} % assumes new font selection scheme installed
\usepackage{amsmath} % assumes amsmath package installed
\usepackage{amssymb}  % assumes amsmath package installed
\usepackage{cite}
\usepackage{subfig}
\usepackage{multirow}
\usepackage{enumitem}
\newtheorem{proposition}{Proposition}[section]

\newcommand*{\Z}{\mathbb{Z}}

\title{Protection Placement for State Estimation Measurement Data Integrity}

\author{Kin Cheong Sou \thanks{The author is with the Department of Electrical Engineering at National Sun Yat-sen University, Taiwan. {\tt sou12@mail.nsysu.edu.tw}. This work is supported by the Ministry of Science and Technology (MOST) of Taiwan through the project: Application of graph decomposition theory in the optimal analysis and planning of large-scale smart electricity grids, MOST 106-2218-E-110-004-MY3.}}

\begin{document}
\maketitle

\begin{abstract}
In this paper a protection placement problem for power network measurement system data integrity is considered. The placement problem is motivated from the data secure power network design applications in [Dan and Sandberg 2010]. The problem is shown to be NP-hard and an integer linear programming formulation is provided based on the topological observability condition by [Krumpholz, Clements and Davis 1980]. The incorporation of the observability condition requires a graph connectivity constraint which is found to be most effectively described by the Miller-Tucker-Zemlin (MTZ) conditions. For a specialization without line power flow measurements, the protection placement problem can be modeled by a domination type integer linear program much easier to solve than the general formulation requiring graph connectivity. Numerical studies with IEEE benchmark and other large power systems with more than 2000 buses indicate that using the proposed formulations, the protection placement problem can be solved in negligible amount of time in realistic application settings.
\end{abstract}

\begin{IEEEkeywords}
Power system state estimation, false data injection attack, equipment placement, observability, integer linear programming. 
\end{IEEEkeywords}

\section{Introduction}
State estimation is a crucial functionality in energy management system. Data integrity of state estimation is of paramount importance to other downstream functionalities including optimal power flow, unit commitment and contingency analysis. Customarily bad data detection (BDD) is employed in conjunction with state estimation to detect and isolate possible data anomalies and to rectify the data if possible. In traditional BDD, however, data anomalies are typically treated as random measurement error or switch status error. This assumption on anomalies significantly limits the energy management system's capabilities to combat malicious data attack on state estimation. For example, references \cite{LRN09,dan2010stealth,sandberg2010security,bobba2010detecting,kosut2010malicious,6504815_TSG,kim2011strategic} analyze and quantify the consequence of a type of malicious data attack called false data injection attack \cite{LRN09}. Another important data integrity research direction, which is the focus of this paper, is protection placement. A typical placement objective is to seek a minimum cost placement of protection resources (e.g., encryption devices, secure phasor measurement units (PMU)) so that, according to the considered attack and defense model, no false data injection attack is possible (e.g., \cite{dan2010stealth, bobba2010detecting, kim2011strategic}).  Because of the combinatorial feature, it is often deemed acceptable to only sub-optimally solve the protection placement problem. For example, \cite{dan2010stealth, bobba2010detecting, kim2011strategic} consider various types of heuristic algorithms aiming to minimize the protection cost. Reference \cite{6504815_TSG} provides a suboptimal strategy to guard against some given attack scenarios. The protection placement problem is closely related to the classical problem of observable measurement system design in power systems (e.g., \cite{aminifar2010contingency, chakrabarti2009placement}). Leveraging the connection between the protection placement problem and observable measurement system design problem, this paper seeks to solve the former problem exactly via integer programming with the help of the results from the later problem (e.g., \cite{KCD80, aminifar2010contingency}). Unifying the results in \cite{KCD80, aminifar2010contingency}, this paper formulates a protection placement optimization problem suitable for applications with bus power injection, line power flow as well as PMU measurements. It is shown that the optimization problem is NP-hard and an integer program formulation of the problem is provided. A crucial consideration in the formulation is the modeling of a graph connectivity requirement due to \cite{KCD80}, which significantly affects the computation effort required to solve the formulation. After extensive investigations it is concluded that the Miller-Tucker-Zemlin (MTZ) conditions \cite{DL91} are best for the purpose herein. The numerical studies in this paper demonstrate that even for models with more than 2000 buses the MTZ based formulation can be solved, on average, within a few minutes on a PC. This is acceptable for power system planning purposes. In addition, for a restricted setup without line power flow measurements (e.g., only PMUs and zero injection buses as in \cite{aminifar2010contingency}) it is established that a streamlined integer linear programming formulation originally appears in \cite{aminifar2010contingency} is indeed correct for the protection placement problem (as well as the problem in \cite{aminifar2010contingency}). To the best of our knowledge, the proof of correctness of the formulation from \cite{aminifar2010contingency} appears for the first time in this paper. Further, by analyzing the structure of the proposed integer programming formulations, we discover that certain binary decision variables can be relaxed to continuous-valued without changing the optimization outcome. Some computation consequences of this discovery are discussed in this paper.

{\bf Outline:} Section~\ref{sec:problem} describes the protection placement problem called perfect protection problem. It is established that the perfect protection problem can be described as an optimization problem with a graph connectivity constraint. In Section~\ref{sec:IP} the perfect protection problem is formulated as a linear integer program with MTZ connectivity constraints. It is also demonstrated that certain binary decision variables can be relaxed to continuous-valued. In Section~\ref{sec:domination} a special setup without line power flow measurement is considered. It is established that the linear integer program derived in Section~\ref{sec:IP} can be specialized to the formulation originally appeared in \cite{aminifar2010contingency}. Section~\ref{sec:studies} describes numerical studies for solving the proposed formulations for realistically-sized examples including the IEEE power system benchmarks. The studies demonstrate the superiority of the proposed formulations. %Section~\ref{sec:conclusion} concludes the paper.

\section{Problem statement} \label{sec:problem}

A power network is modeled as an undirected connected graph where the nodes are buses, and the edges are transmission (distribution) lines. The set of all buses (nodes) is denoted $V$ and the set of lines (edges) is denoted $E$. An edge is defined by an unordered pair $\{i,j\}$ where $i, j \in V$ are the two end nodes of the edge. The power network graph is denoted $(V,E)$. The symbols $|V|$ and $|E|$ denote the number of buses and the number of lines, respectively. Following \cite{dan2010stealth, 6504815_TSG,bobba2010detecting, kim2011strategic}, this paper considers the DC power flow measurement model \cite{Abur_Exposito_SEbook} for state estimation. In this model, the states to be estimated are the bus voltage phasors. The measurements include bus real power injections, line real power flows and phasors measured directly from relevant PMUs. In this paper, a PMU measures the phasor at the installed bus as well as those at all neighboring buses.  Let $z$ and $\theta$ denote the vectors of measurements and states respectively. Then, DC power flow model specifies $z = H \theta + \Delta z$, where $H$ is the measurement matrix
%For example, for the measurement system in Fig.~\ref{fig:ex_3bus} with unit line impedances, 
%\begin{figure}[h]
%\begin{center}
%\includegraphics[width=60mm]{ex_3bus}
%\caption{A 3 bus system with one injection measurement at bus 1, two line power flow measurements (line $\{1,2\}$ and line $\{1,3\}$) and one PMU at bus 3. Each injection or line power flow measurement corresponds to one row in the measurement matrix $H$. However, a PMU typically corresponds to multiple rows of $H$ (i.e., the phasors at its own and neighboring buses).}
%\label{fig:ex_3bus}
%\end{center}
%\end{figure}
%the corresponding $H$ matrix is 
%$$
%H = \begin{bmatrix} \begin{array}{rrr} 2 & -1 & -1 \\ -1 & 1 & 0 \\ -1 & 0 & 1 \\ 1 & 0 & 0 \\ 0 & 1 & 0 \\ 0 & 0 & 1 \end{array} \end{bmatrix}.
%$$
and $\Delta z$ denotes measurement imperfection. In this paper, $\Delta z$ is the vector of data attacks on measurements.

In state estimation, BDD attempts to detect possible data attack in the measurements (i.e., $\Delta z$). In a typical residual-based BDD scheme, the measurement residual $r$ and the data attack $\Delta z$ are related by $r = (I - H (H^T R^{-1} H)^{-1} H^T R^{-1}) \Delta z$ where $R$ is a given diagonal positive definite matrix. The relation specifies that $r$ is the projection of $\Delta z$ on the left null space of $H$. In other words, if $\Delta z = H \tilde{\theta}$ for some $\tilde{\theta}$ then $r = 0$, resulting in detection evasion. Malicious attack of this form is referred to as false data injection attack \cite{LRN09}. As a countermeasure, \cite{dan2010stealth, 6504815_TSG,bobba2010detecting, kim2011strategic} consider the notion of \emph{perfect protection} where selected measurements or PMUs are protected to prevent false data injection attacks. The protection rule is as follows: injection and line power flow measurements can be protected individually: protecting measurement $z_i$ means $\Delta z_i = 0$. However, protection rule for PMU is different: suppose a PMU is associated with measurements $z_j$ with $j \in \mathcal{J}$ for some row index set of $H$. Protecting the PMU means $\Delta z_j = 0$ for all $j \in \mathcal{J}$. Arranging the protected measurements first, a false data injection attack is of the form
\begin{equation} \label{eqn:protected_dz}
\Delta z = \begin{bmatrix} 0 \\ \Delta z_{\bar{p}} \end{bmatrix} = H \tilde{\theta} = \begin{bmatrix} H_p \\ H_{\bar{p}} \end{bmatrix} \tilde{\theta}, \quad \text{for some $\tilde{\theta}$},
\end{equation}
where $H_p$ and $H_{\bar{p}}$ are the submatrices of $H$ corresponding to the protected and unprotected measurements respectively. In (\ref{eqn:protected_dz}) $H_p \tilde{\theta} = 0$ is due to protection directly. To ensure perfect protection (i.e., $\Delta z = 0$), it is required that 
\begin{equation} \label{eqn:perfect_protection_original}
\text{for any $\tilde{\theta}$,} \;\; H_p \tilde{\theta} = 0 \implies \Delta z_{\bar{p}} = H_{\bar{p}} \tilde{\theta} = 0.
\end{equation}
Condition (\ref{eqn:perfect_protection_original}) is equivalent to the fact that the null space of $H_p$ is a subspace of the null space of $H_{\bar{p}}$. However, neither this condition or (\ref{eqn:perfect_protection_original}) is convenient to impose in an optimization problem (for protection placement). Furthermore, a perfect protection condition dependent on the unprotected measurements (i.e., $H_{\bar{p}}$) is undesirable because the constitution of the unprotected measurements may change during operations. This might result in loss of perfect protection. Consequently, in references including \cite{bobba2010detecting, kim2011strategic} the following sufficient condition for perfect protection is considered
%\emph{Perhaps a discussion for the case without PMU. Discussion on $P \neq \emptyset$}.
\begin{equation} \label{eqn:perfect_protection_AO}
H_p \tilde{\theta} = 0 \implies \tilde{\theta} = 0, \; \text{i.e., $H_p$ has full column rank}.
\end{equation}
Condition (\ref{eqn:perfect_protection_AO}) is algebraic observability condition \cite{KCD80} for the reduced measurement system utilizing only the protected measurements characterized by (the rows of) $H_p$. However, the algebraic observability condition is sensitive to the values of line impedance parameters. Thus, instead of (\ref{eqn:perfect_protection_AO}) this paper considers the \emph{topological observability} condition \cite{KCD80} to guarantee perfect protection. Generically, matrix $H_p$ in (\ref{eqn:perfect_protection_AO}) is
\begin{equation} \label{eqn:Hp}
H_p = \begin{bmatrix} M_I A D A^T \\ M_L D A^T \\ M_P \end{bmatrix},
\end{equation}
where $A \in \mathbb{R}^{|V| \times |E|}$ is the (signed) incidence matrix describing the topology of the power network. The directions of the directed arcs in $A$ are designated arbitrarily. Matrix $D \in \mathbb{R}^{|E|^2}$ is diagonal with nonzero entries being the reciprocals of the line reactances. $M_I$, $M_L$ and $M_P$ are submatrices containing rows of identity matrices of appropriate dimensions, in order to select the protected measurements for injection, line power flows and PMU respectively. Then, condition~(\ref{eqn:perfect_protection_AO}) amounts to $H_p$ in (\ref{eqn:Hp}) having full column rank, for a specific $D$. On the other hand, the topological observability condition considered in this paper requires that
\begin{equation} \label{eqn:perfect_protection_TO}
\begin{array}{l}
H_p(\tilde{D}) := \begin{bmatrix} M_I A \tilde{D} A^T \\ M_L \tilde{D} A^T \\ M_P \end{bmatrix} \; \text{has full column rank $(=|V|)$}, \vspace{1mm}
\\ \text{for almost all diagonal $\tilde{D}.$}
\end{array}
\end{equation}
Condition~(\ref{eqn:perfect_protection_TO}) depends on the topology of the power system graph $(V,E)$ and the locations of the protected measurements. However, it is independent of line parameters. Additionally, it can be shown that (\ref{eqn:perfect_protection_TO}) is equivalent to a graph connectivity condition specified by the following statement:
\begin{proposition} \label{thm:perfect_protection}
Let $G = (V,E)$ be a (power network) graph such that $V = \{v_1,v_2,\ldots,v_{|V|}\}$. Let $A \in \mathbb{R}^{|V| \times |E|}$ be the (signed) incidence matrix of $G$ with arc direction arbitrarily designated. Let $I \subseteq V$, $L \subseteq E$ and $P \subseteq V$ be given (i.e., sets of protected injections, line power flows and PMU phasors respectively). Define the following objects:
\begin{itemize}
\item $M_I$ is a submatrix of a $|V| \times |V|$ identity matrix containing the rows corresponding to $I$.

\item $M_L$ is a submatrix of a $|E| \times |E|$ identity matrix containing the rows corresponding to $L$.

\item $M_P$ is a submatrix of a $|V| \times |V|$ identity matrix containing the rows corresponding to $P$.

\item $E_{P0} := \{\{i, 0\} \mid i \in P\}$.
\end{itemize}
Then, condition~(\ref{eqn:perfect_protection_TO}) holds if and only if there exists a function $g: I \mapsto E$, with constraint that $i \in g(i)$ for all $i \in I$, such that the graph $(V \cup \{0\}, L \cup E_{P0} \cup g(I))$ is connected, with $g(I) := \{g(i) \mid i \in I\}$. 
\end{proposition}
\begin{IEEEproof}
For any $\tilde{\theta}$, in (\ref{eqn:perfect_protection_TO}) $H_p(\tilde{D}) \tilde{\theta} = 0$ if and only if
\begin{equation} \label{eqn:augmented_AO}
\begin{bmatrix}
M_I A \tilde{D} A^T \!\! + \!\! M_I M_P^T M_P \!\! & -M_I M_P^T M_P {\bf 1} \\ M_L \tilde{D} A^T & 0 \\ M_P & -M_P {\bf 1}
\end{bmatrix} \!\!\! \begin{bmatrix} \tilde{\theta} \\ \theta_0 \end{bmatrix} = 0, \;
\theta_0 = 0,
\end{equation}
%Let $m_I$ and $m_P$ be the number of rows of matrices $M_I$ and $M_P$ respectively. We specialize matrix $B \in \R^{m_I \times m_P}$ as
%$$
%B(k,j) = \left\{ \begin{array}{cl} 1 & \text{if $M_P(j, i(k)) = 1$} \vspace{1mm} \\ 0 & \text{otherwise} \end{array} \right., \quad k = 1,\ldots,m_I,
%$$
%where $i(k)$ is the (only) column index of $M_I$ for the nonzero entry of $M_I$ at row $k$.
%In (\ref{eqn:augmented_AO}), the symbol ${\bf 1}$ denotes the all-one vector. 
where ${\bf 1}$ is column vector of ones of commensurate dimension. Intuitively, $(M_I A \tilde{D} A^T \!\! + \!\! M_I M_P^T M_P) \tilde{\theta} -M_I M_P^T M_P {\bf 1} \theta_0$ is the vector of injections at $I$ in an augmented graph, denoted $G_a$, with an additional node 0 and additional edges in $E_{P0}$. The equation in (\ref{eqn:augmented_AO}) can be rewritten as
\begin{equation} \label{eqn:augmented_AO1}
\bar{H}_p \begin{bmatrix} \tilde{\theta} \\ \theta_0 \end{bmatrix} := \begin{bmatrix} \bar{M}_I \bar{A} \bar{D} \bar{A}^T \\ \bar{M}_L \bar{D} \bar{A}^T \end{bmatrix} \begin{bmatrix} \tilde{\theta} \\ \theta_0 \end{bmatrix} = 0, \quad \theta_0 = 0
\end{equation}
where $\bar{A}$, $\bar{D}$, $\bar{M}_I$ and $\bar{M}_L$ are the analogies in $G_a$ for $A$, $D$, $M_I$ and $M_L$ in $G$. In particular,
\begin{equation} \label{eqn:bar_matrices}
\begin{array}{l}
\bar{A} = \begin{bmatrix} A & M_P^T \\ 0 & -{\bf 1}^T M_P^T \end{bmatrix}, \;\bar{D} = \begin{bmatrix} \tilde{D} & 0 \\ 0 & I \end{bmatrix}, \; \bar{M}_I = \begin{bmatrix} M_I & 0 \end{bmatrix}, \vspace{1mm} \\
\bar{M}_L = \begin{bmatrix} M_L & 0 \\ 0 & I \end{bmatrix}.
\end{array}
\end{equation}
As a result,
$$
\begin{array}{rcl}
\text{(\ref{eqn:perfect_protection_TO})} & \iff & H_p(\tilde{D}) \tilde{\theta} = 0 \implies \tilde{\theta} = 0, \;\; \text{almost all diag $\tilde{D}$} \vspace{1mm} \\
& \iff & (\ref{eqn:augmented_AO1}) \implies \tilde{\theta} = 0, \;\; \text{almost all diag $\tilde{D}$}
\end{array}
$$
Due to (\ref{eqn:bar_matrices}), $\bar{H}_p$ in (\ref{eqn:augmented_AO1}) can be interpreted as the measurement matrix of a power system corresponding to $G_a$ with injection measurements in $I$ and line power flow measurements in $L \cup E_{P0}$. However, $G_a$ does not have any PMU measurement. With this interpretation, the condition that $(\ref{eqn:augmented_AO1}) \implies \tilde{\theta} = 0$ for almost all $\tilde{D}$ is topological observability of the measurement system described with $I$ and $L \cup E_{P0}$ for $G_a$. According to \cite{KCD80}, this condition is equivalent to the existence of a function $g: I \mapsto E$, with constraint that $i \in g(i)$ for all $i \in I$, such that the graph $(V \cup \{0\}, L \cup E_{P0} \cup g(I))$ is connected.
\end{IEEEproof}

Note that since $A^T {\bf 1} = 0$ and node 0 can only be connected through $E_{P0}$, the two equivalent conditions established in Proposition~\ref{thm:perfect_protection} both requires $P \neq \emptyset$ (i.e., at least one protected PMU). However, Proposition~\ref{thm:perfect_protection} does not assume connectivity of the original power network graph $(V,E)$, though connectivity typically holds in practice. In the situation without PMU, condition~(\ref{eqn:perfect_protection_AO}) should be modified so that the null space of $H_p$ is the span of ${\bf 1}$. The topological observability condition in \cite{KCD80} can be utilized directly to describe the modified condition without Proposition~\ref{thm:perfect_protection}.

The protection placement problem in this paper seeks a minimum cost protection placement to guarantee perfect protection in the sense of (\ref{eqn:perfect_protection_TO}). The protection decisions are described by 0-1 binary decision variables $x_i$ for $i \in V$, $y_j$ for $j \in E$ and $z_k$ for $k \in V$ for bus power injections, line power flows and PMU respectively. By convention, $x_i = 1$ if and only if the injection at bus $i$ is protected. Analogous conventions apply to the $y_j$ and $z_k$ decision variables. In addition, let $\mathcal{M}^I \subseteq V$, $\mathcal{M}^L \subseteq E$ and $\mathcal{M}^P \subseteq V$ denote the set of buses with measured injections, the set of lines with measured power flows and the set of buses with PMU installed, respectively. Then, by convention $x_i = 0$ if $i \notin \mathcal{M}^I$. Similarly, $y_j = 0$ if $j \notin \mathcal{M}^L$ and $z_k = 0$ if $k \notin \mathcal{M}^P$. Let $c_i^I$, $c_j^L$ and $c_k^P$ denote the protection costs associated with $x_i$, $y_j$ and $z_k$ respectively. If bus $v$ is a zero-injection bus the corresponding protection cost is $c_v^I = 0$. Then, the objective of the protection placement problem is to minimize
\begin{equation} \label{eqn:obj}
\sum\limits_{i \in \mathcal{M}^I} c_i^I x_i + \sum\limits_{j \in \mathcal{M}^L} c_j^L y_j + \sum\limits_{k \in \mathcal{M}^P} c_k^P z_k.
\end{equation}
Let $I(x) = \{i \in V \mid x_i = 1\}$, $L(y) = \{e \in E \mid y_e = 1\}$ and $P(z) = \{i \in V \mid z_i = 1\; \text{or} \; \exists k \in V : z_k = 1, \; \{i,k\} \in E\}$ be the set of buses with power injection measurement protected, the set of lines with power flow measurement protected and the set of buses whose phasors are measured by protected PMU respectively. Recall $E_{P0}(z) := \{\{i,0\} \mid i \in P(z)\}$. Then, by Proposition~\ref{thm:perfect_protection} the perfect protection constraint is
\begin{equation} \label{eqn:constraint}
\begin{array}{l}
\text{there exists $g: I(x) \mapsto E$ with $i \in g(i)$ for all $i \in I(x)$} \\
\text{s.t.~$(V \cup \{0\}, L(y) \cup E_{P0}(z) \cup g(I(x)))$ is connected.}
\end{array}
\end{equation}
In sequel, (\ref{eqn:obj}) and (\ref{eqn:constraint}) together are referred to as the perfect protection problem. The perfect protection problem is NP-hard. This can be established by considering the special case with a PMU at each bus but nothing else (i.e., $\mathcal{M}^I = \emptyset$, $\mathcal{M}^L = \emptyset$, $\mathcal{M}^P = V$). This special case reduces to a minimum dominating set problem which is NP-hard (e.g., \cite{garey2002computers}).

\section{Integer programming formulation of perfect protection problem} \label{sec:IP}
To model (\ref{eqn:obj}) and (\ref{eqn:constraint}) as an integer program it is necessary to describe the connectivity requirement in (\ref{eqn:constraint}) as linear constraints with respect to the decision variables. Intuitively, (\ref{eqn:constraint}) is equivalent to the existence of a rooted spanning tree using edges in $L(y) \cup E_{P0}(z) \cup g(I(x)))$. First, for an arbitrary graph $(\bar{V},\bar{E})$ the conditions for existence of a rooted spanning tree are described. Second, additional requirements are described specifying that the spanning tree can be formed using only edges in $L(y) \cup E_{P0}(z) \cup g(I(x)))$, making the spanning tree conditions relevant to (\ref{eqn:constraint}).

For existence of spanning tree well-known conditions include, for example, Miller-Tucker-Zemlin (MTZ) conditions \cite{MTZ}, subtour elimination (e.g., \cite{LLKS85}), Martin's conditions \cite{Martin91}, single commodity flow (e.g., \cite{LFRR12, SL17}) and multi-commodity flow (e.g., \cite{NW12}). After computation studies conducted in conjunction with this work (not shown in this paper), a variant of MTZ conditions \cite{DL91} is adopted. The detail is as follows. Let $(\bar{V}, \bar{E})$ be given and $\bar{V}$ includes some ``reference'' node 0. Let $\bar{\mathcal{A}}$ be the bi-directed version of $\bar{E}$ (i.e., $\bar{\mathcal{A}} := \{(i,j), (j,i) \mid \{i,j\} \in \bar{E}\}$). Let $f_{ij} \in \{0,1\}$ for $(i,j) \in \bar{\mathcal{A}}$ be 0-1 binary decision variables (to describe the spanning tree). In addition, define integer decision variables $u_i \in \{0,1,\ldots,|\bar{V}|-1\}$ for $i \in \bar{V}$ (the $u$ variables can be relaxed to real numbers between 0 and $|\bar{V}|-1$ \cite{pataki2003teaching}). Then, according to \cite{DL91} the graph $(\bar{V}, \bar{E})$ contains a (directed) spanning tree if and only if there admits $f$ and $u$ satisfying the following constraints
\begin{subequations}
\begin{align}
& (|\bar{V}| \! - \!1) f_{ij} \! + \!(|\bar{V}| \!-\!3) f_{ji} \!+ \!u_i \!- \!u_j \le (|\bar{V}|\!-\!2), \; \forall (i,j) \in \bar{\mathcal{A}} \label{subeqn:ste} \\
& 0 \le u_i \le |\bar{V}| - 1, \quad \forall i \in \bar{V} \label{subeqn:label} \\
& \sum\limits_{j: (i,j) \in \bar{\mathcal{A}}} f_{ij} = 1, \quad \forall i \in \bar{V} \setminus \{0\} \label{subeqn:outgoing} \\
& \sum\limits_{j: (j,0) \in \bar{\mathcal{A}}} f_{j0} \ge 1 \label{subeqn:incoming} \\
& \sum\limits_{(i,j) \in \bar{\mathcal{A}}} f_{ij} = |\bar{V}| - 1 \label{subeqn:totalE}
\end{align}
\end{subequations}
The directed spanning tree, if exists, is characterized by $f$ (i.e., $f_{ij} = 1$ if and only if arc $(i,j)$, from $i$ to $j$, is included). The constraints (\ref{subeqn:ste}) to (\ref{subeqn:totalE}) impose restrictions in the choice of $f$ (and $u$). Constrants~(\ref{subeqn:ste}) and (\ref{subeqn:label}) together are referred to as subtour elimination constraints. They prevent directed cycles with less than $|\bar{V}|$ arcs. Constraint~(\ref{subeqn:outgoing}) specifies that for each node other than 0 exactly one outgoing arc is included in $f$. Constraint~(\ref{subeqn:incoming}) specifies that for node 0 at least one incoming arc is included in $f$. Node 0 is the in-root of the spanning tree. Constraint~(\ref{subeqn:totalE}) requires that $f$ contains $|\bar{V}|$ arcs. It can be verified that $f$ corresponding to any directed spanning tree rooted at node 0 satisfies all constraints in (\ref{subeqn:ste}) to (\ref{subeqn:totalE}), with some appropriate choice of $u$. Conversely, $f$ (together with $u$) satisfying constraints (\ref{subeqn:ste}) to (\ref{subeqn:incoming}) means that from each node other than node 0 there exists a directed path to node 0. Hence, $f$ characterizes a connected graph containing all nodes. Constraint (\ref{subeqn:totalE}) states that the connected graph contains exactly $|\bar{V}|-1$ arcs which implies that it is a tree (cf.~\cite{diestel2005graph}). In conclusion, constraints (\ref{subeqn:ste}) to (\ref{subeqn:totalE}) are necessary and sufficient for existence of (directed) spanning tree in $(\bar{V}, \bar{E})$ (or equivalently the connectedness of $(\bar{V}, \bar{E})$).

Next, we specialize conditions (\ref{subeqn:ste}) to (\ref{subeqn:totalE}) to $\bar{V} = V \cup \{0\}$ and $\bar{E} = E \cup \{\{i,0\} \mid i \in V\}$. In addition, we impose the conditions that $f_{ij}$ is enabled (i.e., allowed to be one) if and only if $\{i,j\} \in L(y) \cup E_{P0}(z) \cup g(I(x)))$ for some function $g$ in (\ref{eqn:constraint}). An edge of the form $\{i, 0\}$ for $i \in V$ can only be enabled by $E_{P0}(z)$. Recall that PMU at bus $i$ is protected if and only if $z_i = 1$. Therefore,
\begin{equation}
f_{0i} + f_{i0} \le z_i + \sum\limits_{j : \{i,j\} \in E} z_j, \quad \forall i \in V.
\end{equation}
On the other hand, an edge $e \in E$ is enabled either because $e \in L(y)$ (i.e., protected line) or $e \in g(I(x))$ (i.e., assignment from protected injection):
\begin{equation} \label{eqn:E_original_enabling}
f_{ij} + f_{ji} \le y_{\{i,j\}} + w_{\{i,j\}i} + w_{\{i,j\}j}, \quad \forall \{i,j\} \in E,
\end{equation}
where $w_{ei}$ for $e \in E$, $i \in V$ are 0-1 binary decision variables describing the assignment $g(I(x))$ in (\ref{eqn:constraint}). The convention is that $w_{ei} = 1$ if and only if $g(i) = e$. Because of the fact that $g$ is a function, as well as other requirements in (\ref{eqn:constraint}), additional constraints should be imposed on $w$:
\begin{equation} \label{eqn:assignment}
w_{ei} = 0, \quad i \notin e, \quad \text{and} \quad \sum\limits_{e \in E} w_{ei} \le x_i, \quad \forall i \in V.
\end{equation}
This concludes modeling (\ref{eqn:constraint}) as linear constraints. The integer program describing the perfect protection problem is summarized as
\begin{equation} \label{opt:perfect_protection}
\begin{array}{cl}
\underset{x, y, z, w, f, u}{\text{minimize}} & \text{cost in (\ref{eqn:obj})} \vspace{2mm} \\
\text{subject to} & \text{constraints (\ref{subeqn:ste}) to (\ref{eqn:assignment})} \vspace{2mm} \\
& x_i = 0, \; \forall i \notin \mathcal{M}^I, \;\; x_i \in \{0,1\}, \; \forall i \in \mathcal{M}^I \vspace{2mm} \\
& y_j = 0, \; \forall j \notin \mathcal{M}^L, \;\; y_j \in \{0,1\}, \; \forall j \in \mathcal{M}^L \vspace{2mm} \\
& z_k = 0, \; \forall k \notin \mathcal{M}^P, \;\; z_k \in \{0,1\}, \; \forall j \in \mathcal{M}^P \vspace{2mm} \\
& w_{ei} \in \{0,1\}, \quad \forall (e,i) \in E \times V, \vspace{2mm} \\
& f_{ij}, f_{ji} \in \{0,1\}, \; \forall \{i,j\} \in E \cup \{\{k, 0 \} \! \mid \! k \in V\} \vspace{2mm} \\
& u_i \in \{0,\ldots,|V|\}, \quad \forall i \in \{0\} \cup V
\end{array}
\end{equation}
In constraints (\ref{subeqn:ste}) to (\ref{subeqn:totalE}) in problem~(\ref{opt:perfect_protection}) the convention is that $\bar{V} = V \cup \{0\}$ and $\bar{\mathcal{A}} = \{(i,j), (j,i) \mid \{i,j\} \in \bar{E}\}$ with $\bar{E} = E \cup \{\{i,0\} \mid i \in V\}$. 

Some of the integer variables in (\ref{opt:perfect_protection}) can be relaxed to continuous variables. As noted in \cite{pataki2003teaching}, the $u$ variables can be continuous, between 0 and $|V|$. In addition, the coefficient matrices multiplying the $x$ variables and $y$ variables are submatrices of identity. By a total unimodularity argument (e.g., \cite{schrijver1998theory}) it is seen that $x$ and $y$ can be relaxed to continuous variables between 0 and 1. Furthermore, the $w$ variables appear only in constraints~(\ref{eqn:E_original_enabling}) and (\ref{eqn:assignment}). With respect to $w$ these are constraints describing an assignment problem (e.g., \cite{schrijver1998theory}). Thus, the $w$ variables can also be relaxed to continuous ones between 0 and 1. In summary, the integer program in (\ref{opt:perfect_protection}) is equivalent to the following mixed integer linear program
\begin{equation} \label{opt:perfect_protection_mixed}
\begin{array}{cl}
\underset{x, y, z, w, f, u}{\text{minimize}} & \text{cost in (\ref{eqn:obj})} \vspace{2mm} \\
\text{subject to} & \text{constraints (\ref{subeqn:ste}) to (\ref{eqn:assignment})} \vspace{2mm} \\
& x_i = 0, \; \forall i \notin \mathcal{M}^I, \;\; 0 \le x_i \le 1, \; \forall i \in \mathcal{M}^I \vspace{2mm} \\
& y_j = 0, \; \forall j \notin \mathcal{M}^L, \;\; 0 \le y_j \le 1 \; \forall j \in \mathcal{M}^L \vspace{2mm} \\
& z_k = 0, \; \forall k \notin \mathcal{M}^P, \;\; z_k \in \{0,1\}, \; \forall j \in \mathcal{M}^P \vspace{2mm} \\
& 0 \le w_{ei} \le 1, \quad \forall (e,i) \in E \times V, \vspace{2mm} \\
& f_{ij}, f_{ji} \in \{0,1\}, \; \forall \{i,j\} \in E \cup \{\{k, 0 \} \! \mid \! k \in V\} \vspace{2mm} \\
& 0 \le u_i \le |V|, \quad \forall i \in \{0\} \cup V
\end{array}
\end{equation}

\section{Domination specialization of perfect protection problem} \label{sec:domination}
While (\ref{opt:perfect_protection}) describes the perfect protection problem in general, more streamlined (hence less computationally demanding) formulations are available if restrictions are imposed on the available measurement types (characterized by $\mathcal{M}^I$, $\mathcal{M}^L$ and $\mathcal{M}^P$). An example mentioned earlier is that $\mathcal{M}^I = \emptyset$, $\mathcal{M}^L = \emptyset$ and $\mathcal{M}^P = V$. In this case,  problem~(\ref{opt:perfect_protection}) is a minimum dominating set problem which requires fewer constraints to model and less time solve. It turns out a more general case with $\mathcal{M}^L = \emptyset$ also admits streamlined formulation as a domination-type problem without requiring graph connectivity. This is specified by the following statement
\begin{proposition} \label{thm:no_line}
Let graph $(V,E)$ be given with $V = \{v_1, \ldots, v_{|V|}\}$. Let $I \subseteq V$ and $P \subseteq V$ be given (i.e., sets of protected injections and protected PMU phasors respectively). Consider the following two conditions:
\begin{enumerate}[label=(\alph*)]
\item there exists $g : I \mapsto E$ such that $i \in g(i)$ and the graph $(V \cup \{0\}, g(I) \cup E_{P0})$ is connected, with $E_{P0} = \{\{i,0\} \mid i \in P\}$,
\item there exists $h : I \mapsto V$ satisfying two requirements: (i) $h(i) = i$ or $\{i, h(i)\} \in E$ and (ii) $P \cup h(I) = V$.
\end{enumerate}
Then, if (a) holds then (b) holds. Conversely, if $(V,E)$ is connected, $P \neq \emptyset$ and (b) holds then (a) holds.
\end{proposition}
\begin{IEEEproof}
Suppose (a) holds. Connectivity of graph $(V \cup \{0\}, g(I) \cup E_{P0})$ implies that there exists a spanning tree rooted in node 0. Let $T$ denote the tree. We will show that $T$ can be used to define an node-to-node assignment function $h$ so that (b) holds. Initialize $h(v) := \emptyset$ for all $v \in I$. Consider a reverse level order traversal of $T$ and then truncate the traversal by removing all leaves of $T$ and node 0. Let $i_1, i_2, \ldots, i_m$ be the resulted truncation. Consider inspecting $i = i_1, i_2, \ldots, i_m$ one-at-a-time. For any node $i$ inspected, let $C(i) \subseteq V$ be the set of children of $i$ in $T$ ($C(i) \neq \emptyset$ since $i$ is not a leaf). The nodes in $\{i\} \cup C(i)$ induce a subgraph of $T$ with $|C(i)|$ edges of the form $\{i,j\}$ for $j \in C(i)$. Since $i \neq 0$, (a) in the statement implies that these edges must be in $g(\{i\} \cup C(i))$. In addition, for each $j \in C(i)$ either (I) $j \in I$ and $g(j) = \{i,j\}$ or (II) $j \notin I$ or $g(j) \neq \{i,j\}$ (then $i \in I$ and $g(i) = \{i,j\}$). For case (I), update $h(j) := j$. For case (II), update $h(i) := j$. This finishes the inspection for node $i$. The process is repeated for all nodes in the truncated tree traversal. We will prove the following three claims: (A) for all $v \in I$, $h(v)$ is updated at most once, (B) for all $v \in I$, either $h(v) = v$ or $\{v, h(v)\} \in E$ and (C) for all $u \in V \setminus P$, there exists some $i_u \in I$ such that $h(i_u) = u$. Claim (A) is due to the update rule of $h$: updating $h(v) := w$ for some $w \in V$ means that $g(v) = \{v,w\}$. By (a), $g(v)$ corresponds to only one edge. In addition, the edge $\{v,w\}$ is considered exactly once during the inspection process. Claim (B) is a direct consequence of the update rule for $h$. To see claim (C), note that when $u$ is considered (when its parent is inspected), $i_u$ is identified and the corresponding $h(i_u)$ is updated to $u$. By claim (A), assignment $h(i_u) = u$ will not be modified again during subsequent inspections. Hence, claim (C) holds. Claims (B) and (C) together imply that $h$ satisfies (i) and (ii) of (b) in the statement. This shows that (a) implies (b).

Suppose (b) holds, $(V,E)$ is connected and $P \neq \emptyset$. To show (a) we will show that it is possible to construct a node-to-edge assignment function $g$ such that the graph $(V \cup \{0\}, g(I) \cup E_{P0})$ in (a) contains a path from any $i \in V$ to node 0. Initialize $R_0 := P \cup \{0\}$ (the set of nodes ``reaching'' node 0 and $g(v) := \emptyset$ for all $v \in I$. Next, consider $i \notin R_0 \supseteq P$, if any. By (b), there exists some $j \in I$ such that $\{i,j\} \in E$ and $h(j) = i$. However, $j$ need not be unique. Pick an arbitrary $j$ and denote $j := h^{-1}(i)$. If $j \in R_0$ update $g(j) := \{i,j\}$ and update $R_0 \leftarrow R_0 \cup \{i\}$. This is motivated by the fact that now a path exists from $i$ to $R_0$ (and hence to node 0), with edges in $g(I) \cup E_{P0}$ . On the other hand, if $j \notin R_0$ then again by (b) there exists some $k \in I$ such that $\{j,k\} \in E$ and $k = h^{-1}(j) := h^{-2}(i)$. Continue to identify $h^{-3}(i), \ldots$. Since graph $(V,E)$ is finite, at some point this process must stop because of two possibilities (P1): $h^{-d}(i) \in R_0$ for some $d \in \Z$ or (P2): $h^{-d}(i) = i$ for some $d$. Note that the following variant of (P2) is impossible: $h^{-d}(i) = h^{-s}(i) = v \neq i$ for some $s \ge 1$. In this case, $h(v)$ would have two different values, violating (b) in the statement. For possibility (P1), denote $h^0(i) = i$ and update $g(h^{-\ell}(i)) := \{h^{-\ell+1}(i), h^{-\ell}(i)\}$ for $\ell = 1,\ldots,d$. Update also $R_0 \leftarrow R_0 \cup \{i, h^{-1}(i), \ldots, h^{-d+1}(i)\}$. By updating $g$, the nodes $i, h^{-1}(i), \ldots, h^{-d+1}(i)$ become connected to $h^{-d}(i) \in R_0$ with edges in $g(I) \cup E_{P0}$ and hence they are included in $R_0$. On the other hand, possibility (P2) means that $i$, $h^{-1}(i)$, \ldots, $h^{-d}(i) = i$ form a cycle with every node in $I$. This concludes the process with $i \notin R_0$. Next, find another node not in (the most updated) $R_0$ or the cycle(s) identified so far, and repeat the same process of updating $g$ and $R_0$ until $V$ is partitioned into $R_0$ and possibly a finite number of cycles. The assumption $P \neq \emptyset$ implies $R_0 \neq \emptyset$. In addition, connectivity of $(V,E)$ implies that, if cycles exist, at least one is adjacent to $R_0$ in the sense that there exist node $u \in R_0$ and node $v$ in the cycle such that $\{u,v\} \in E$. Suppose the cycle is $v = h^0(v), h^{-1}(v), \ldots, h^{-d}(v) = v$. Then update $g(h^{-\ell}(v)) := \{h^{-\ell}(v), h^{-\ell+1}(v)\}$ for $\ell = 1,\ldots,d-1$ and $g(h^{-d}(v)) := \{u,v\}$. Update $R_0$ to include all the nodes in the cycle and continue with the next cycle adjacent to the updated $R_0$ until $R_0 = V \cup \{0\}$ (otherwise adjacent cycle must exist). We will prove the following three claims: (D) for all $v \in I$, $g(v)$ is updated at most once, (E) for all $v \in I$, $v \in g(v)$ and (F) while $R_0$ and $g$ are being updated, at any moment $(R_0, g(I) \cup E_{P0})$ is connected. To see claim (D), if $v$ is to be included in $R_0$ then $g(v)$ is updated for the first time because a node is included in $R_0$ at most once. Only in possibility (P1) can $v$ be already in $R_0$ while $g(v)$ is being updated. However, for (P1) $g(v) := \{v, u\}$ for some $u$ implies $u$ is to be included in $R_0$ and $h(v) = u$. If $g(v)$ is updated more than once, for example, as $g(v) := \{v,u\}$ and $g(v) := \{v,w\}$, then $u \neq w$ (since a node is included in $R_0$ at most once) and $h(v)$ would be valued $u$ and $w$ simultaneously. This violates (b). Thus, claim (D) is true. Claim (E) is a direct consequence of the update rule for $g$. For claim (F), note that initially when $R_0 = P \cup \{0\}$, graph $(R_0, g(I) \cup E_{P0})$ is connected because of $E_{P0}$. When some $i, h^{-1}(i), \ldots, h^{-d}(i)$ is included in $R_0$, these nodes are connected to $R_0$ through edges $g(h^{-\ell}(i))$ for $\ell = 1,\ldots,d$. In addition, by claim (D) the values of $g(h^{-\ell}(i))$ for $\ell = 1,\ldots,d$ will not change in subsequent updates of $g$. Hence, claim (F) is true. Since $R_0 = V \cup \{0\}$ at the end, claims (E) and (F) imply that $g$ satisfies condition (a) in the statement. In conclusion, assuming $P \neq \emptyset$ and $(V,E)$ connected, (b) implies (a) in the statement.
\end{IEEEproof}
Since it is impossible to satisfy (a) in Proposition~\ref{thm:no_line} if $P = \emptyset$, the perfect protection problem can be augmented with the constraint $P \neq \emptyset$. In addition, assume that $(V,E)$ is connected (connectivity can be checked in polynomial time). Then, conditions (a) and (b) in Proposition~\ref{thm:no_line} are equivalent. Thus, (b) provides an alternative formulation of the perfect protect problem in the case with no line power flows (i.e., $\mathcal{M}^L = \emptyset$). The protection decisions are encoded by 0-1 binary variables $x$ and $z$, for injections and PMUs respectively, as in (\ref{eqn:obj}) and (\ref{eqn:constraint}). The sets $I(x)$ and $P(z)$ are defined similarly. In condition (b), $I = I(x)$ and $P = P(z)$. To describe function $h$, define 0-1 binary decision variables $w_{ij}$ for $i,j \in V$ such that $w_{ij} = 1$ if and only if $h(j) = i$. The conditions $P(z) \cup h(I(x)) = V$ and $h(i) = i$ or $\{i, h(i)\} \in E$ require
\begin{equation} \label{eqn:domination_z}
z_i + \sum\limits_{k : \{i,k\} \in E} z_k + \sum_{j \in V} w_{ij} \ge 1, \quad \forall i \in V.
\end{equation}
as well as
\begin{equation} \label{eqn:domination_w}
\sum\limits_{i \in V} w_{ij} \le x_j, \; \forall j \in V, \quad w_{ij} = 0, \; \forall i \neq j : \{i,j\} \notin E.
\end{equation}
While constraint $P \neq \emptyset$ is a natural consequence of (a) in Proposition~\ref{thm:no_line} (i.e., (\ref{eqn:constraint})), it needs to be explicitly imposed while formulating with (b) in Proposition~\ref{thm:no_line}:
\begin{equation} \label{eqn:P_not_empty}
\sum\limits_{k \in \mathcal{M}^P} z_k \ge 1.
\end{equation}
Then, the following integer program is equivalent to problem (\ref{eqn:obj}) and (\ref{eqn:constraint}) when $\mathcal{M}^P = \emptyset$ and $(V,E)$ is connected:
\begin{equation} \label{opt:perfect_protection_domination}
\begin{array}{cl}
\underset{x, z, w}{\text{minimize}} & \sum\limits_{i \in \mathcal{M}^I} c_i^I x_i + \sum\limits_{k \in \mathcal{M}^P} c_k^P z_k \vspace{2mm} \\
\text{subject to} & \text{constraints~(\ref{eqn:domination_z}), (\ref{eqn:domination_w}) and (\ref{eqn:P_not_empty})} \vspace{2mm} \\
& x_i = 0, \; \forall i \notin \mathcal{M}^I, \;\; x_i \in \{0,1\}, \; \forall i \in \mathcal{M}^I \vspace{2mm} \\
& z_k = 0, \; \forall k \notin \mathcal{M}^P, \;\; z_k \in \{0,1\}, \; \forall j \in \mathcal{M}^P \vspace{2mm} \\
& w_{ij} \in \{0,1\}, \; \forall i,j \in V
\end{array}.
\end{equation}
Constraint~(\ref{eqn:domination_z}) requires that each node is dominated by at least a protected PMU or a protected injection (subject to assignment rule), motivating the term ``domination'' for problem (\ref{opt:perfect_protection_domination}). Empirically, the graph connectivity requirement in (\ref{opt:perfect_protection}) is much more difficult to handle than the domination requirement in (\ref{opt:perfect_protection_domination}). The difference in computation performances will be demonstrated in  Section~\ref{sec:studies}. Note that by the same argument from problem~(\ref{opt:perfect_protection}) to (\ref{opt:perfect_protection_mixed}), the $x$ and $w$ variables in (\ref{opt:perfect_protection_domination}) can be relaxed to continuous variables between 0 and 1 resulting in the following equivalent mixed integer linear program
\begin{equation} \label{opt:perfect_protection_domination_mixed}
\begin{array}{cl}
\underset{x, z, w}{\text{minimize}} & \sum\limits_{i \in \mathcal{M}^I} c_i^I x_i + \sum\limits_{k \in \mathcal{M}^P} c_k^P z_k \vspace{2mm} \\
\text{subject to} & \text{constraints~(\ref{eqn:domination_z}), (\ref{eqn:domination_w}) and (\ref{eqn:P_not_empty})} \vspace{2mm} \\
& x_i = 0, \; \forall i \notin \mathcal{M}^I, \;\; 0 \le x_i \le 1, \; \forall i \in \mathcal{M}^I \vspace{2mm} \\
& z_k = 0, \; \forall k \notin \mathcal{M}^P, \;\; z_k \in \{0,1\}, \; \forall j \in \mathcal{M}^P \vspace{2mm} \\
& 0 \le w_{ij} \le 1, \; \forall i,j \in V
\end{array}.
\end{equation}

In \cite{aminifar2010contingency} problem (\ref{opt:perfect_protection_domination}) is used to model a PMU placement problem equivalent to (\ref{eqn:obj}) and (\ref{eqn:constraint}). However, no proof is provided in \cite{aminifar2010contingency} to justify the use. To the best of our knowledge, this paper is the first in the literature to establish that (\ref{opt:perfect_protection_domination}) indeed models the perfect protection problem (as well as the one in \cite{aminifar2010contingency}), under the relevant assumptions discussed earlier.

\begin{figure}[h]
\begin{center}
\includegraphics[width=60mm]{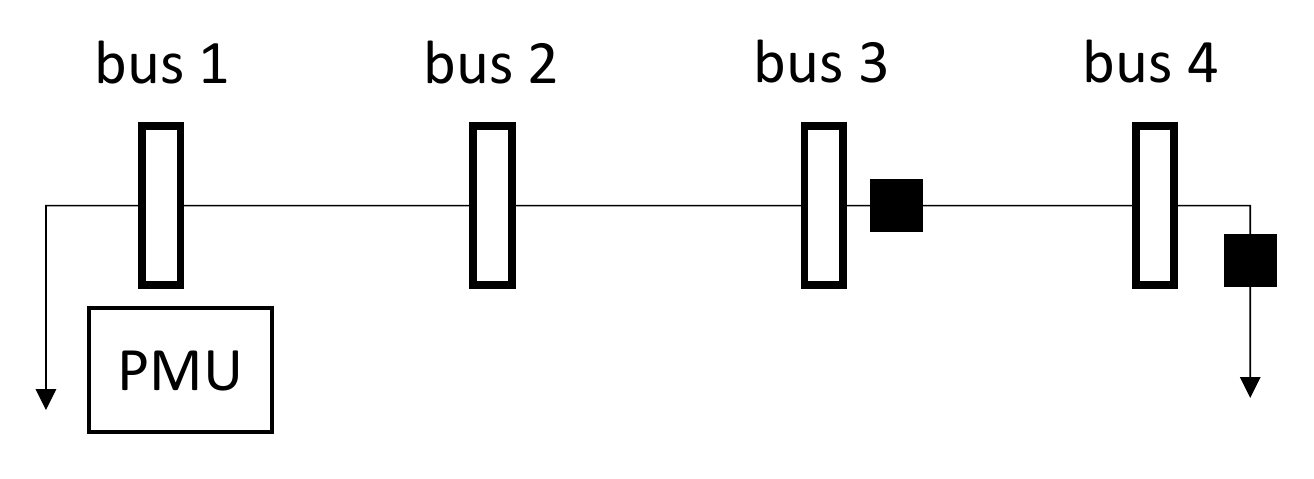}
\caption{All three measurements/device are protected but perfect protection is not achieved. The line power flow and injection are in fact the same measurement, and hence (\ref{eqn:perfect_protection_TO}) cannot be satisfied. However, per the conjecture every node is dominated (e.g., 1, 2 by PMU, 3 by line power flow and 4 by injection).}
\label{fig:counterex_Aminifar}
\end{center}
\end{figure}

One might conjecture that the following extension of (\ref{opt:perfect_protection_domination}) may model the perfect protection problem in general. In addition to domination due to nodes with protected PMU or protected injection, a function $f : E \mapsto V$ such that $f(e) \in e$ can assign an edge with protected line power flow to one of its end nodes to let it dominated. However, the conjecture is not true. See Fig.~\ref{fig:counterex_Aminifar} for a counterexample.

% \section{Cuts}

% \section{Minimum spanning tree}

\section{Numerical studies} \label{sec:studies}

\begin{figure}[!h]
\begin{center}
\includegraphics[width=80mm]{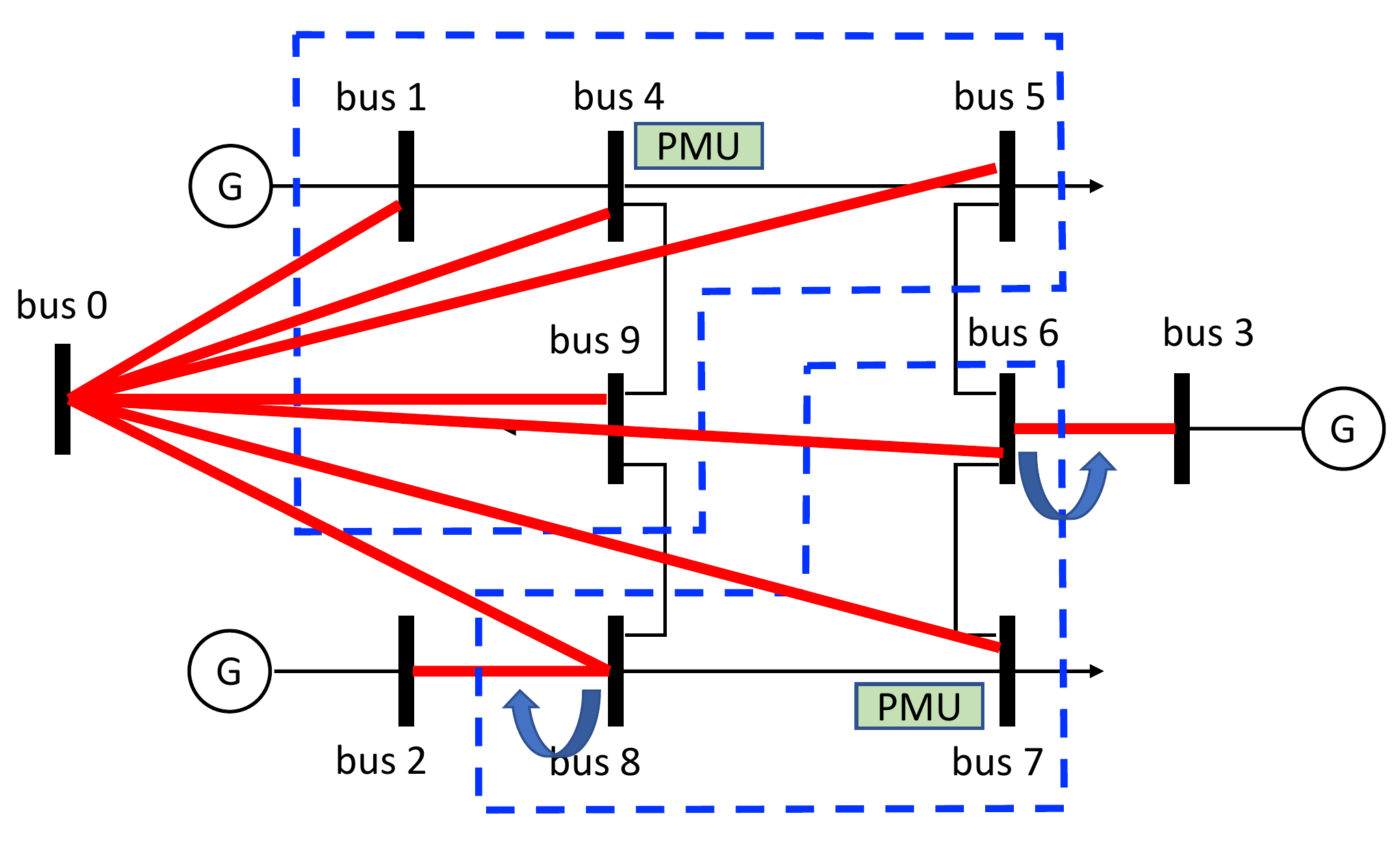}
\caption{Minimum cost perfect protection layout of the IEEE 9 bus example, in which the PMUs at bus 4 and bus 7 are protected. The result is obtained by solving problem~(\ref{opt:perfect_protection}). In the figure, the two dash blue boxes indicate the ``zones'' covered by the protected PMUs at bus 4 and bus 7 respectively. The two zones enable direct connections to the ``reference'' bus 0 through the thick red lines for buses 1, 4, 5, 9, 6, 7 and 8. In addition, buses 4, 6 and 8 are zero-injection buses whose injections are known to be zero (i.e., protected automatically). These injections can help to satisfy (\ref{eqn:constraint}) in the form of $g(I)$ (i.e., being assigned to an incident line). For example, the protected zero-injection at bus 6 is assigned to line $\{3,6\}$ (thick red) and the one at bus 8 is assigned to line $\{2,8\}$ (thick red). In summary, constraint (\ref{eqn:constraint}) is satisfied by the red thick lines being the edges in $E_{P0}(z) \cup g(I(x))$.}
\label{fig:IEEE9_Krumpholz}
\end{center}
\end{figure}

\begin{figure}[!h]
\begin{center}
\includegraphics[width=80mm]{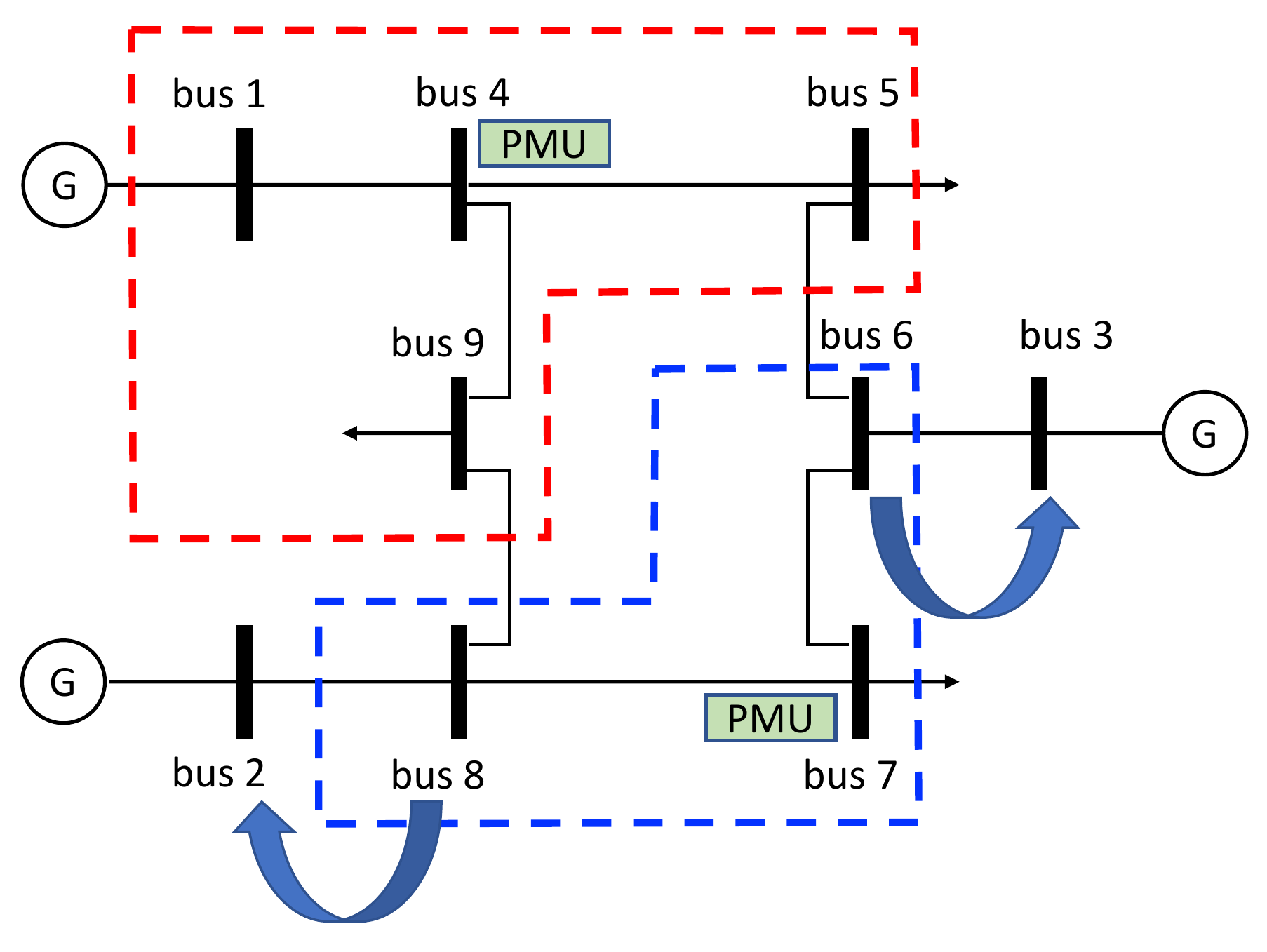}
\caption{Minimum cost perfect protection layout of the IEEE 9 bus example, in which the PMUs at bus 4 and bus 7 are protected. The result is obtained by solving problem~(\ref{opt:perfect_protection_domination}). In the figure, the blue and red dash boxes indicate the ``zones'' dominated by protected PMUs at bus 4 and bus 7 respectively, in the sense that the buses in the zones are in set $P$ as in (b) of Proposition~\ref{thm:no_line}. In addition, the two zero-injection buses can be ''assigned'' through function $h$ in (b) of Proposition~\ref{thm:no_line} to dominate the two buses not in any zone (i.e., buses 2 and 3).}
\label{fig:IEEE9_Aminifar}
\end{center}
\end{figure}

To illustrate the perfect protection layout and the computation experience with the proposed formulations, instances of problems~(\ref{opt:perfect_protection}), (\ref{opt:perfect_protection_mixed}), (\ref{opt:perfect_protection_domination}) and (\ref{opt:perfect_protection_domination_mixed}) are solved in this section. All optimization problems are solved using Gurobi \cite{gurobi} in MATLAB and the computations are performed on a PC with 14 CPU cores at 2GHz with 128GB of RAM. The power network graphs describing the instances are from the IEEE power system benchmark database \cite{IEEE_power_grid_benchmarks}. First the IEEE 9-bus system is considered. In this example, all bus injections, line power flows and the PMUs at all buses are eligible for protection. The protection costs for all injections and line flows are one except at the three zero-injection buses (i.e., 4, 6, 8) the injection protection cost is zero. On the other hand, the PMU protection cost is one for all buses. The protection costs are specified so that at optimality only the injections at the three zero-injection buses and the PMUs would be candidates for protection. Because of Proposition~\ref{thm:no_line} and the choice of protection costs, the perfect protection problem in this example is equivalent to the minimum cost PMU placement problem in \cite{aminifar2010contingency}. Two formulations in (\ref{opt:perfect_protection}) and  (\ref{opt:perfect_protection_domination}) are solved for the perfect protection problem in this example. Note that (\ref{opt:perfect_protection_domination}) is applicable since the line power flows will never be protected at optimality. Both cases lead to the same result: protecting the PMUs at bus 4 and bus 7 with the total protection cost being two. This is the same cost for the same example in \cite{aminifar2010contingency}, though the result therein is to ``protect'' the PMUs at bus 5 and bus 8. It can be verified that, with identity $\tilde{D}$, matrix $H_p(\tilde{D})$ in (\ref{eqn:perfect_protection_TO}) has full column rank with protected PMUs at bus 4 and bus 7 (with three zero-injection buses at 4, 6 and 8).
%$$
%H_p = \begin{bmatrix} \begin{array}{rrrrrrrrr}
%    -1  &   0  &   0  &   3  &  -1 &    0 &    0  &   0 &   -1 \\
%     0   &  0  &  -1  &   0  &  -1 &    3  &  -1  &   0 &    0 \\
%     0  &  -1  &   0  &   0  &   0  &   0  &  -1  &   3 &  -1 \\
%     1   &  0  &   0  &  0   &  0   &  0  &   0   & 0   &  0 \\
%     0  &   0  &   0  &   1  &   0  &   0  &   0  &   0 &    0 \\
%     0  &   0  &   0  &   0  &   1  &   0  &   0  &   0  &   0 \\
%     0  &   0  &   0  &   0  &   0  &   1  &   0  &   0  &   0 \\
%     0  &   0  &   0  &   0  &   0  &   0  &   1  &   0  &   0 \\
%     0  &   0  &   0  &   0  &   0  &   0  &   0  &   1  &   0 \\
%     0  &   0  &   0  &   0  &   0  &   0  &   0  &   0  &   1
%     \end{array}
%\end{bmatrix}_{10 \times 9}.
%$$
Fig.~\ref{fig:IEEE9_Krumpholz} shows the optimal protection layout and verifies that constraint~(\ref{eqn:constraint}) (i.e., the only constraint of perfect protection problem) is satisfied.
 Fig.~\ref{fig:IEEE9_Aminifar} shows how condition (b) in Proposition~\ref{thm:no_line} is satisfied by the optimal protection layout due to solving (\ref{opt:perfect_protection_domination}).

\begin{figure}[!t]
\begin{center}
\includegraphics[width=80mm]{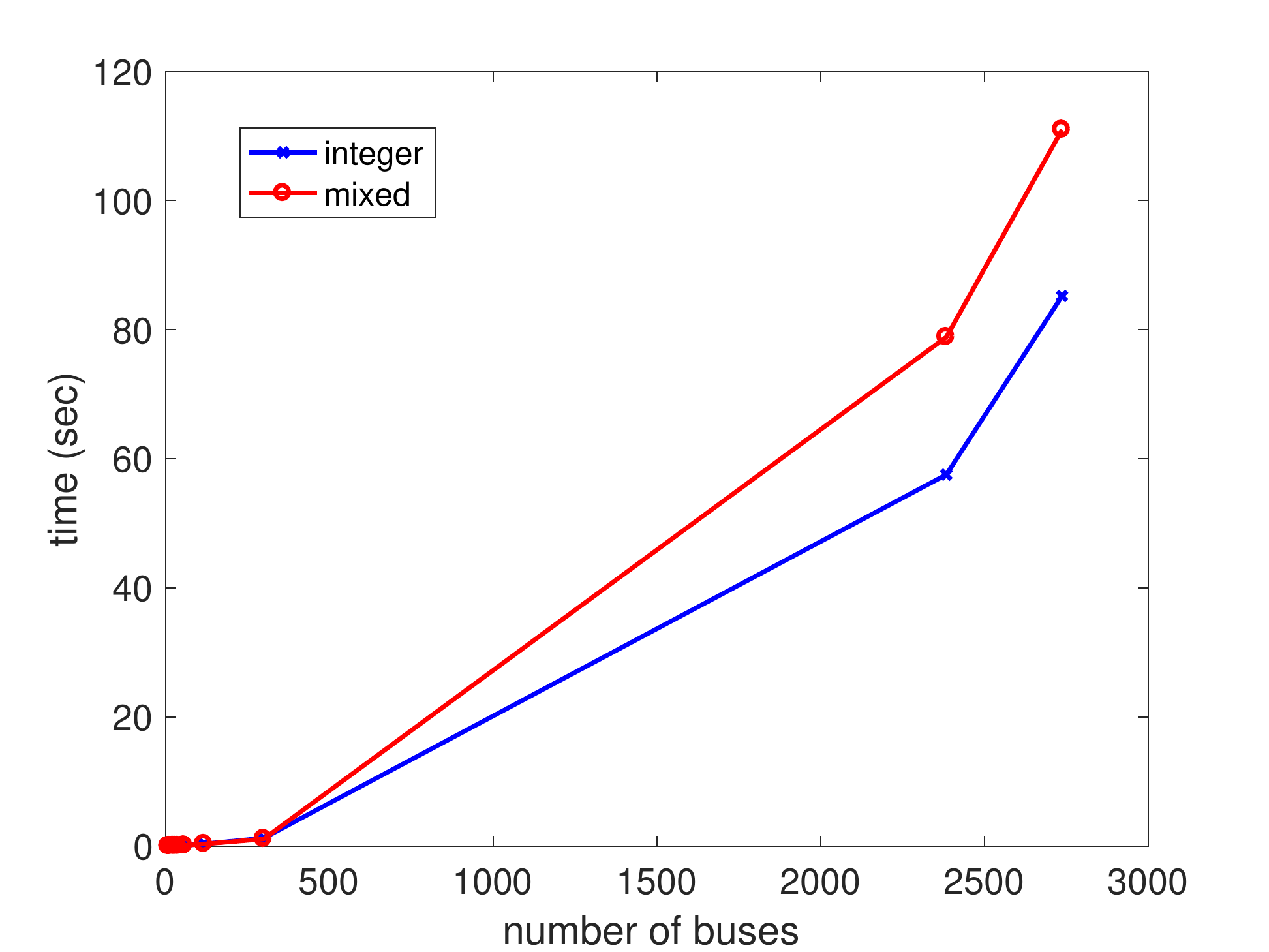}
\caption{Average time (over 100 samples) for solving (\ref{opt:perfect_protection}) (i.e., integer program) and (\ref{opt:perfect_protection_mixed}) (i.e., mixed integer program) for IEEE power system benchmarks of various numbers of buses: 9, 14, 24, 30, 39, 57, 118, 300, 2383 and 2736.}
\label{fig:full}
\end{center}
\end{figure}

The remaining study is divided into two parts referred to as (i) the full measurement case and (ii) the no line power flow case pertaining Section~\ref{sec:IP} and Section~\ref{sec:domination} respectively. In the full measurement case, all bus injections and line power flows are measured. In addition, each bus is equipped with a PMU and 10\% of the buses are  zero-injection buses. That is, $\mathcal{M}^I = V$, $\mathcal{M}^L = E$ and $\mathcal{M}^P = V$ and there exists $\mathcal{Z} \subseteq V$ such that $c_i^I = 0$ for $i \in \mathcal{Z}$ and $|\mathcal{Z}| \approx 0.1 |V|$. The injection protection cost $c_i^I$ for $i \in \mathcal{M}^I \setminus \mathcal{Z}$ and the line protection cost $c_j^L$ for line $j \in \mathcal{M}^L$ are random integers uniformly sampled between 1 and 100. The PMU protection cost $c_k^P$ for bus $k \in \mathcal{M}^P$ is 80\% of the sum of the line protection costs of all incident lines, rounded up to the next integer. The setup for the no line power flow case is similar except that the line power flow measurements are not protected (i.e., $\mathcal{M}^L = \emptyset$). For the full measurement case, for each power system in database \cite{IEEE_power_grid_benchmarks}, 100 instances of the general perfect protection problem are generated and the corresponding integer program in (\ref{opt:perfect_protection}) and mixed integer program in (\ref{opt:perfect_protection_mixed}) are solved. For all instances, solving the two formulations results in the same (optimal) protection cost. In addition, the protection layouts are verified to provide perfect protection as specified by (\ref{eqn:perfect_protection_TO}), for some randomly chosen $\tilde{D}$. For each power system, the average computation time (over 100 samples) for solving (\ref{opt:perfect_protection}) and (\ref{opt:perfect_protection_mixed}) is shown in Fig.~\ref{fig:full}. All instances, including 2383 bus and 2736 bus cases, can be solved within two minutes on average in the experiment. Further, it is noted that for larger systems the mixed integer formulation (\ref{opt:perfect_protection_mixed}) in fact requires more time to solve than the pure integer formulation (\ref{opt:perfect_protection}).
However, for smaller systems the opposite is true. This is illustrated in Table~\ref{tab:full_small}.
\begin{table}[h]
\caption{Average computation time (sec) for solving integer program in (\ref{opt:perfect_protection}) and mixed integer program in (\ref{opt:perfect_protection_mixed})} \label{tab:full_small}
\begin{tabular}{|c|c|c|c|c|c|c|}
\hline
$|V|$ & 9 & 14 & 30 & 57 & 118 & 300 \\
\hline
integer (\ref{opt:perfect_protection}) & 0.0105 &   0.0215 &  0.0426 & 0.1182 &   0.3534  &  1.2446 \\
\hline
mixed (\ref{opt:perfect_protection_mixed}) & 0.0105  &  0.0224 & 0.0419 & 0.1122  &  0.3433 & 1.1143 \\
\hline
\end{tabular}
\end{table}
%\begin{figure}[h]
%\begin{center}
%\includegraphics[width=80mm]{full_zoomed.eps}
%\caption{Zoomed in version of Fig.~\ref{fig:full} showing more clearly the comparison of the average computation time for smaller instances below 300 buses. The mixed integer formulation in (\ref{opt:perfect_protection_mixed}) requires less time to solve.}
%\label{fig:full_zoomed}
%\end{center}
%\end{figure}

\begin{figure}[!t]
\begin{center}
\includegraphics[width=80mm]{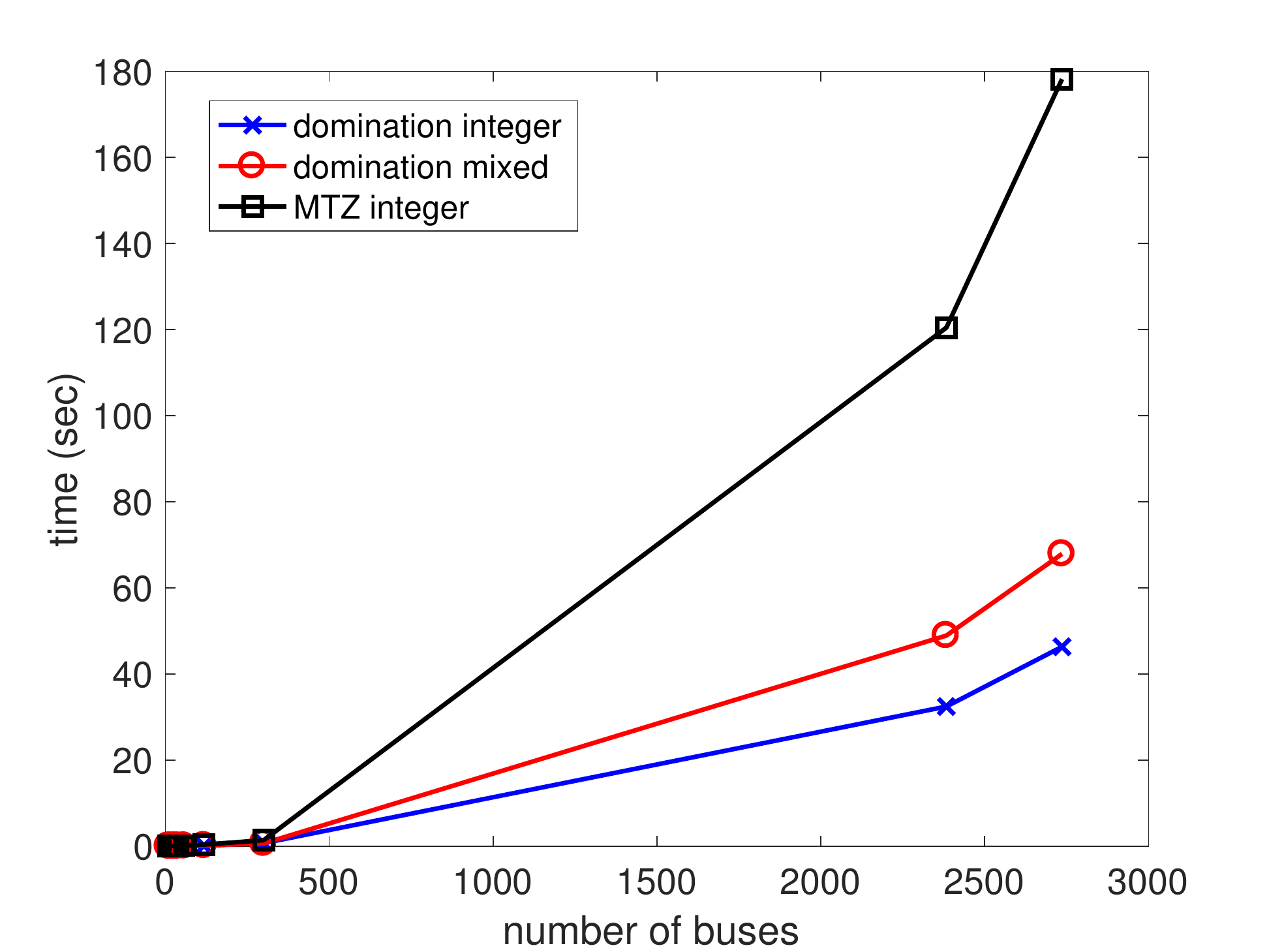}
\caption{Average time (over 100 samples) for solving (\ref{opt:perfect_protection_domination}) (i.e., domination integer program), (\ref{opt:perfect_protection_domination_mixed}) (i.e., domination mixed integer program) and (\ref{opt:perfect_protection}) (i.e., general formulation) for IEEE power system benchmarks of various numbers of buses: 9, 14, 24, 30, 39, 57, 118, 300, 2383 and 2736.}
\label{fig:domination}
\end{center}
\end{figure}

For the no line power flow case 100 random instances of the perfect protection problems are generated for each power system, and these instances are solved as the corresponding domination formulations in (\ref{opt:perfect_protection_domination}) for the pure integer version and (\ref{opt:perfect_protection_domination_mixed}) for the mixed integer version. In addition, the general formulation in (\ref{opt:perfect_protection}) is solved for comparison. For all instances, all three formulations result in the same (optimal) protection cost. The average computation times (over 100 samples) for instances of various numbers of buses are shown in Fig.~\ref{fig:domination}.
The results indicate that the streamlined domination formulations in (\ref{opt:perfect_protection_domination}) and (\ref{opt:perfect_protection_domination_mixed}) (applicable only in the no line power flow case) require only a fraction of the time required by the general formulation in (\ref{opt:perfect_protection}). Further, up to the 300 bus case on average it takes less time to solve the mixed integer formulation in (\ref{opt:perfect_protection_domination_mixed}) (not shown in detail). However, for larger instances the purely integer formulation in (\ref{opt:perfect_protection_domination}) requires less time to solve.
%Fig.~\ref{fig:domination_zoomed}.
%\begin{figure}[h]
%\begin{center}
%\includegraphics[width=80mm]{domination_zoomed.eps}
%\caption{Zoomed in version of Fig.~\ref{fig:domination} showing more clearly the comparison of the average computation time for smaller instances below and including 300 buses. The mixed integer formulation in (\ref{opt:perfect_protection_domination_mixed}) requires less time to solve.}
%\label{fig:domination_zoomed}
%\end{center}
%\end{figure}

\section{Conclusion} \label{sec:conclusion}
Extending the result in \cite{KCD80}, it is possible to model as integer program with graph connectivity constraints the perfect protection (or observable measurement placement) problem including power injections, line power flows as well as PMUs. The MTZ constraints modeling graph connectivity are found to be superior for computation efficiency. Even for instances with more than 2000 buses optimal placement can be found within a few minutes on a PC. This agrees with the computation findings in \cite{NW12} for variants of the dominating set problem. For the case without line power flow measurements the perfect protection problem can be reduced into a domination type integer program, which further improves computation efficiency. As a byproduct, the result in this paper proves the correctness of the formulation in \cite{aminifar2010contingency}. Also, the studied integer programs contain certain binary variables that can be relaxed to continuous variables. However, the computation advantage of the relaxations has not been strongly established by the studies conducted in this paper. A further investigation in this issue could be beneficial. In addition, extension of the results in this paper to general $N-k$ contingency cases could open up more application opportunities in the future.

\bibliographystyle{IEEEtran}
\bibliography{SDIP}
\end{document}